\documentclass[onecolumn]{IEEEtran}
\usepackage{amsmath}

\usepackage{cuted}

\usepackage{array}
\usepackage{textcomp}
\usepackage{nicefrac}
\usepackage{amssymb}
\usepackage{amsfonts}
\usepackage{caption}
\usepackage{multirow}
\usepackage{cite}
\usepackage{textcomp}
\usepackage{graphicx}
\usepackage{subcaption}
\usepackage{epsfig}
\usepackage{epstopdf}
\usepackage{bm}

\usepackage{microtype} 
\usepackage{balance}
\usepackage{stfloats}

\usepackage{float} 
\usepackage{hyperref} 
\usepackage{amsmath}
\usepackage{color}
\usepackage{graphicx}
\usepackage{algorithm,algorithmicx,algpseudocode}
\graphicspath{ {images/} }

\begin{document}
\title{
	Energy Storage Management via Deep Q-Networks 
} 


\author{
	Ahmed S. Zamzam,
	Bo Yang, and
	Nicholas D. Sidiropoulos
	\thanks{The work of A.S. Zamzam and N.D. Sidiropoulos was supported in part by NSF grant 1525194. 
		A. S. Zamzam and B. Yang are with the ECE Dept., Univ. of Minnesota, Minneapolis, MN 55455, USA. N. D. Sidiropoulos is with the Department of Electrical and Computer Engineering, University of Virginia, Charlottesville, VA 22904.
		E-mails: ahmedz@umn.edu, yang4173@umn.edu, nikos@virginia.edu.
	}
}
\maketitle

\begin{abstract}
Energy storage devices represent environmentally friendly candidates to cope with volatile renewable energy generation. Motivated by the increase in privately owned storage systems, this paper studies the problem of real-time control of a storage unit co-located with a renewable energy generator and an inelastic load. Unlike many approaches in the literature, no distributional assumptions are being made on the renewable energy generation or the real-time prices. Building on the deep Q-networks algorithm, a reinforcement learning approach utilizing a neural network is devised where the storage unit operational constraints are respected. The neural network approximates the action-value function which dictates what action (charging, discharging, etc.) to take. Simulations indicate that near-optimal performance can be attained with the proposed learning-based control policy for the storage units. 

\end{abstract}


\section{Introduction}
Renewable energy generation has expanded rapidly in order to enable the environmental, social and economic benefits of future smart grids.
Due to their intermittent nature, renewable energy sources increase the variability of the energy generation portfolio and introduce new challenges of real-time control. Energy storage systems provide an opportunity to increase the grid efficiency and cope with fluctuating renewable energy generators~\cite{Braff-2016}. The ability of batteries to respond quickly to changes in the system makes them ideal candidate for a wide range of power systems applications, including spatio-temporal energy arbitrage~\cite{Pandzi-2015}, peak shaving~\cite{Gayme-2013}, frequency~\cite{Makarov-2012}, and congestion management~\cite{Vargas-2015,Pandzi-2015}. 

This paper focuses on optimal energy management for a single power consumer (or microgrid system) with an energy storage unit, a renewable energy source, and an inelastic load demand. Privately owned and operated storage units are being rapidly installed by both large-scale energy consumers~\cite{guo2011cutting} as well as individuals~\cite{peterson2010economics}. In addition, social benefits have been identified for energy consumer ownership of energy storage units~\cite{sioshansi2010welfare}. 

The operation of energy storage units has recently received significant attention. One line of research concerning this problem assumes a deterministic time-varying environment~\cite{Makarov-2012, lamont2013assessing, Zamzam-2018}. Thus, the operational and planning decisions are made by solving an off-line optimization problem that takes into account operational losses and/or installation costs.

On the other hand, several prior works~\cite{atwa2010optimal, grillo2012optimal, rahbar2015real, xu2017optimal} have studied real-time or online energy management under stochastic assumptions regarding the environment. The authors of~\cite{atwa2010optimal,grillo2012optimal} assumed 
full knowledge of idealized
stationary stochastic processes for the demand and the renewable energy generation. However, even under stationary assumptions, the knowledge of the distribution is very hard to obtain by the consumer. In~\cite{rahbar2015real}, a real-time management framework was proposed that combines an off-line optimal solution with a sliding-window optimization that requires estimates of the demand, renewable energy generation, and prices. 
Estimating these using only consumer-side data is one challenge; complexity is another. An analytic solution was sought in~\cite{xu2017optimal} where it was shown that the optimal solution can be characterized using two thresholds. The consumer charges the energy storage until it is fully charged when the energy price is below a certain threshold, and discharges it until it is empty if the price is above another threshold. Albeit the approach is simple and well-motivated, the two thresholds can be obtained only if the demand and energy prices are known in advance.

Reinforcement learning (RL) is an area of machine learning that deals with an agent learning how to behave in an environment by performing actions and observing the results. Classical RL, e.g., Q-learning~\cite{sutton1998reinforcement}, considers an environment with discrete states which is not suitable for the power systems application where the state of the system is mostly continuous. Combining advances in deep learning~\cite{krizhevsky2012imagenet} with Q-learning resulted in {\it deep Q-networks} (DQN)~\cite{mnih2015human} which achieved human-level control for several Atari video games. In a nutshell, deep neural networks were proposed for estimating the action-value function.

In this paper, we propose a learning-based framework for operating energy storage units building on the DQN algorithm which matches the continuous nature of the energy management system. The proposed approach learns a control policy using interactions with the environment with no distributional knowledge of the load, the renewable energy generation, or the power prices. A neural network is fitted utilizing historical observations to approximate the action-value function. Therefore, the consumer can implement the action with the highest expected reward. The control policy ensures that the battery operational constraints are never violated. The simulation results show near-optimal performance of the proposed approach.

The remainder of this paper is organized as follows. In Section II, the energy management system model is presented. Then, the DQN approach is introduced in Section III. The proposed framework is put forth in Section IV. In Section V, the simulation results are presented, and the paper is concluded in Section VI.

\section{Energy Management System Modeling}
\begin{figure}[h]
	\centering
	\includegraphics[width=0.25\textwidth]{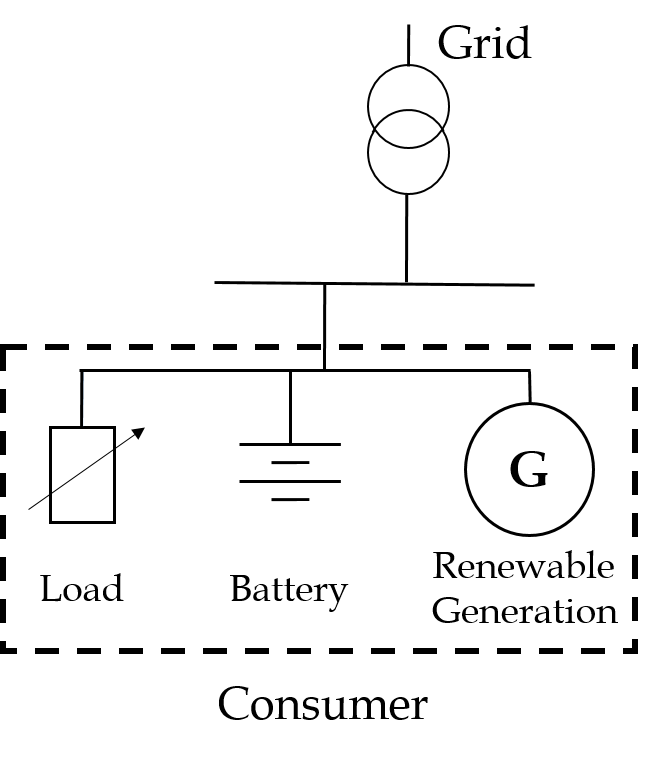}
	\caption{Architecture of the considered model.}
	\label{fig:system}
\end{figure}

This paper studies the operation of a storage unit owned by an electricity consumer. The charging and discharging of the storage unit is assumed to be locally controlled. The consumer purchases power from the grid to charge the battery, or sells its stored energy to the grid. The temporal axis is discretized into time slots with duration $ \delta $. We denote the state of charge of the battery at time slot $ t $ by $ e_t $. In addition, let the capacity of the battery be denoted by $ E $ (in KWh). The consumer also has a renewable energy source installed. Fig.~\ref{fig:system} depicts the model considered in this paper.

For each time slot $ t $, the price of purchasing power from the grid is denoted by $ p_t $ (in \$/KWh). It is assumed that the consumer can also sell energy to the grid at the same price. The inelastic power demand of the consumer at time slot $ t $ is denoted by $ d_t $. Notice that there is no specific stochastic process assumed for driving the energy price or demand, and the two quantities are available to the consumer only at the beginning of time slot $ t $. A renewable energy source is installed at the consumer side where $ r_t $ represents the power generated at time slot $ t $. Again, this quantity is observed by the consumer only at the beginning of time slot $ t $. The consumer can utilize this energy to satisfy the inelastic demand $ d_t $, charge the storage unit, or sell the excess power to the grid. 

During time slot $ t $, the consumer purchases $ b_t $ amount of electricity for the energy storage. If $ b_t $ is negative, then the consumer sells energy to the grid. For the rest of the paper, we denote by $ b_t^+ $ and $ b_t^- $ the positive and negative parts of $ b_t $, respectively. Then, $ b_t^+ $ and $ b_t^- $ can be written as
$$ b_t^+ = \max\{0, b_t\},\qquad b_t^- = - \min\{0, b_t\}. $$
Let $ \overline{r}_c $ and $ \overline{r}_d $ denote the maximum charging and discharging rates, respectively. Therefore, the charging and discharging rates are bounded as follows
\begin{equation}\label{eq:rate-con}
\eta_c \ b_t^+ \leq \overline{r}_c, \qquad \frac{b_t^-}{\eta_d} \leq \overline{r}_d
\end{equation}   
where $ \eta_c\! \in (0, 1]$ and $ \eta_d\!\in (0, 1] $ are the charging and discharging efficiency coefficients. Therefore, the change in the state of charge of the battery can be written as follows
\begin{equation}\label{eq:Dyn_battery}
e_{t+1} = e_t + \delta\ \eta_c\ b_t^+ - \delta\ \frac{b_t^-}{\eta_d}
\end{equation}
where $ \delta $ is the slot duration.
The decision to purchase $ b_t $ power from the grid is {\it feasible} if $ 0 \leq e_{t+1} \leq E $ and~\eqref{eq:rate-con} is satisfied.

The consumer receives a utility $ u_t( b_t + d_t - r_t, p_t ) $ at time slot $ t $, which depends on the net energy withdrawn from the grid $ (b_t + d_t - r_t) $ and the energy price $ p_t $. Then, the energy storage is controlled in order to minimize the discounted cost function given by
\begin{equation}\label{eq:Cost}
C_{t} = \mathbb{E}\Big\{ \sum_{t' = t}^{\infty} \gamma^{t'}\ u_{t'}( b_{t'} + d_{t'} - r_{t'},\ p_{t'} ) \Big\}
\end{equation}
where $ \gamma \in (0, 1) $ denotes the discount factor. Let $ \mathbf{s}_t := [d_t - r_t,\ p_t,\ e_t]^T $ denote the {\it state} of the system at time slot $ t $. Then, we define a {\it policy} $ \boldsymbol{\pi} := [\mu_0, \mu_1, \ldots] $ as a sequence of decision rules such that $ \mu_t({\bf s}_t) $ is feasible for all $ t $. 

\section{Deep Reinforcement Learning}
\begin{figure}[h]
	\centering
	\includegraphics[width=0.45\textwidth]{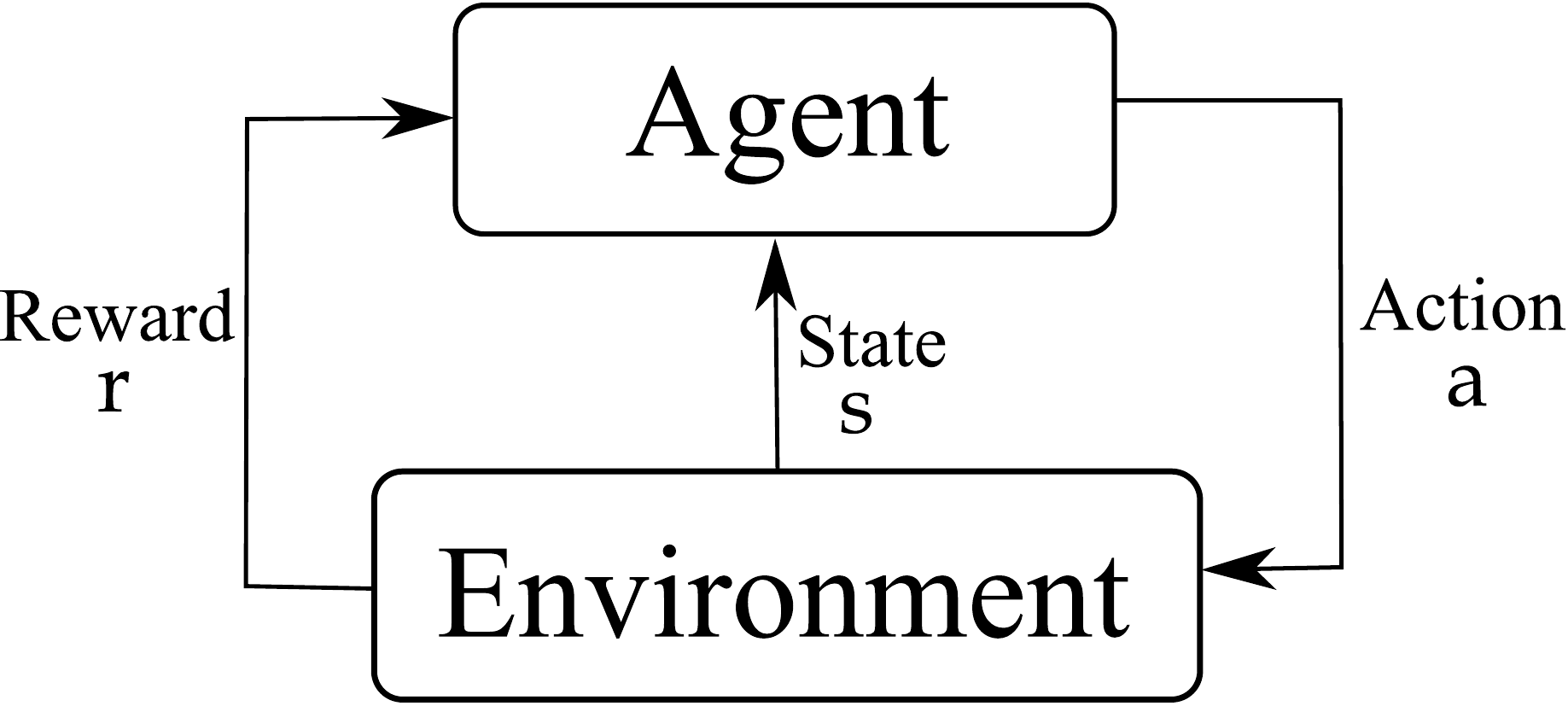}
	\caption{Reinforcement learning process.}
	\label{fig:RL}
\end{figure}

Reinforcement learning is a branch of machine learning that is concerned with optimizing the performance of an agent in an environment; see Fig.~\ref{fig:RL}. The interaction between the agent and the environment is represented using three main elements which are state, actions, and reward. The state $ {\bf s} $ is a vector collecting the state of the exogenous and endogenous environment parameters. The actions are a set of interactions that the agent can perform at any given system state. Finally, the reward is a function that encodes the agent's desired performance and provides a quantifiable immediate response of the environment to the implementation of a specific action in a state.

Classical reinforcement learning methods, e.g., Q-learning, were developed to handle environments with discrete states and actions. Therefore, the action-value of the actions in all the possible states can be described using a look-up table. Upon observing the state of the environment, the agent performs the action that has the highest action-value. However, in many physical systems, the state of the environment is naturally continuous, and hence, the classic Q-learning approaches are not directly applicable.

The deep reinforcement learning approach~\cite{mnih2013playing, mnih2015human} utilizes a neural network that approximates the action-value function $ Q({\bf s}, a) $ which is defined as the maximum expected return achievable by following any strategy. That is, the optimal action-value function can be expressed as
\begin{equation}
Q^\star ({\bf s}, a) = \max_{\boldsymbol{\pi}}\  \mathbb{E}\big[ \sum_{t' = t}^{\infty} \gamma^{t'-t}\ r_{t'} | {\bf s}_t = {\bf s}, a_t = a, \boldsymbol{\pi}\big]. 
\end{equation}
Known as universal function approximators, neural networks are perfectly suited for accurate representation of the action-value function. A deep neural network with weights $ \boldsymbol{\theta} $ is used to approximate $ Q({\bf s}, a) $ which is referred to as a Q-network.

	\begin{algorithm}[t]                
		\caption{Deep Q-Networks Learning}
		\label{alg:DQN}    
		\vspace{0.2cm}
		\begin{minipage}{40cm}
			\indent\textbf{Input:} $0<\gamma<1$, $ 0<\epsilon<0$, and $ 0<\kappa<1$
			\vspace{0.2cm}
			\begin{algorithmic}[1]
				\State\emph{Intialization:} Replay memory $ \mathcal{D} $, and action-value\\ \hspace{2cm}function $ Q({\bf s}, a; \boldsymbol{\theta}) $ with random weights $ \boldsymbol{\theta}_0 $.
				\vspace{0.1cm}
				\State 	\textbf{while} {$ t \leq T $}
				\textbf{do}\hspace{0.1cm}
				\State \hspace{.5cm} With probability $ \epsilon $, select a random action $ a_t $.
				\State \hspace{.5cm} Otherwise, select $ a_t = \arg\max_{a'} Q({\bf s}_t, a'; \boldsymbol{\theta}_t) $.
				\State \hspace{.5cm} Implement action $ a_t $.
				\State \hspace{.5cm} Observe the reward $ r_t $ and the new system state $ {\bf s}_{t+1} $.
				\State \hspace{.5cm} Store the tuple $ ({\bf s}_t, a_t, r_t, {\bf s}_{t+1}) $ in $ \mathcal{D} $.
				\State \hspace{.5cm} Sample a random minibatch $ \mathcal{J} $ from $ \mathcal{D} $.
				\State \hspace{.5cm} Set $ y_j = r_j + \gamma \max_{a'} Q({\bf s}_{j+1}, a'; \boldsymbol{\theta}_t) $ for all $ j \in \mathcal{J} $.
				\State \hspace{.5cm} Perform gradient step on $ \sum_{j\in\mathcal{J}} (y_j - Q({\bf s}_j, a_j; \boldsymbol{\theta}_t))^2 $.
				\State \hspace{.5cm} Update $ \epsilon = \kappa\epsilon $.
				\State\textbf{end while}
			\end{algorithmic}
		\end{minipage}
	\end{algorithm}

The optimal action-value function $ Q^{\star}({\bf s}, a) $ obeys the {\it Bellman equation} which can be written as
\begin{equation}
Q^\star( {\bf s}_t, a_t ) = r_t + \gamma \max_{a'} Q^{\star} ( {\bf s}_{t+1}, a').
\end{equation}
The idea behind the training of reinforcement learning methods is to use Bellman equation to devise an iterative update rule. Such update rule eventually converges to the optimal action-value function~\cite{sutton1998reinforcement} under mild conditions on the step size. The learning procedure is summarized in Algorithm~\ref{alg:DQN} where the $ \epsilon $-greedy approach is used to ensure adequate exploration of the state space.

During the learning process, we utilize the so-called {\it experience replay} technique~\cite{lin1993reinforcement}. This approach stores tuples representing past experiences in a replay memory. Hence, at each update, a randomly drawn minibatch from past experiences are selected to perform gradient updates on the neural network parameterizing the action-value function. Algorithm~\ref{alg:DQN} shows the implementation of the experience replay mechanism within our learning procedure.

\section{Energy Storage Management using DQN}
In this section, we present the learning-based energy storage management system. 
At any time slot, the consumer chooses to charge the battery, discharge it, or to set $ b_t $ to be zero. The consumer utilizes a deep neural network to approximate the action-value function $ Q({\bf s}, a) $.

The consumer aims at minimizing the discounted cost function~\eqref{eq:Cost} using the DQN approach. In the proposed approach, after sufficient exploration of the state space, the consumer implements the action $ a' $ that maximizes $ Q({\bf s}, a) $. Afterwards the consumer gets a reward $ r_t $ for the implemented action, and observes the new state of the system $ {\bf s}_{t+1} $. Then, using Bellman equation the value of the $ Q({\bf s}_t, a') $ at optimality should be 
$ r_t + \gamma \max_a Q({\bf s}_{t+1}, a) $.
Therefore, a training sample is created to update the parameters of the neural network using the observed reward and the new state of the system as in Algorithm~\ref{alg:DQN}.

The proposed learning-approach is {\it model-free} in the sense that no assumptions are being made on the processes generating the energy prices or the renewable energy generation. The approach learns the mapping using samples obtained by interactions with the system. 

Utilizing the DQN approach, the control policy $ \boldsymbol{\pi} $ of the storage unit is defined such that 
\begin{equation}
a_t = \mu_t({\bf s}_t) = \arg\max_{a'} Q({\bf s}_t,\ a') 
\end{equation}
where the function $ Q({\bf s}_t, a) $ is defined by the deep neural network. In order to ensure feasibility of the actions implemented using this policy, the rate of charge and discharge have to be chosen such that the battery is not over charged nor over discharged.
When the selected action is to charge the storage unit, the rate of  charge is set to be the maximum charging rate $ \overline{r}_c $ unless the state of charge is expected to exceed the maximum storage capacity if battery is charged with this rate for a time duration of length $ \delta $. Therefore, the charging rate is set to be
$$ b_t^+  =  \max \{\ 0,\ \min \{\ \overline{r}_c,\ (E - e_t)\ /\ \delta \eta_c  \}$$
which ensures that the state of charge of the battery is always kept below $ E $. Similarly, when the selected action is to discharge the battery, the discharge rate is set to be
$$ b_t^-  =  \min \{\ \overline{r}_d,\ \max \{\ 0,\ e_t \eta_d\ /\ \delta  \} \}$$
which ensures that there is enough energy available to discharge.

\section{Experiments}
We now present the results of our experiments to assess the performance of the learning-based energy management framework. We use a time slot duration $ \delta = 5 $ minutes. A measured load data is utilized to simulate the customer consumption of energy. The load data was measured at $ 18 $ buses, with a $ 1 $-second resolution, from a feeder in Anatolia, California, during a day in the summer of 2012~\cite{bank2013development}. In order to get enough training data, we concatenate the load profiles at the $ 18 $ buses in the feeder which results in a sequence of power consumption of $ 18 $ days with resolution of $ 1 $ second. Therefore, for each period of $ 5 $ minutes duration, we have $ 300 $ load measurements. Then, a total of $ 25 $ sequences with a resolution of $ 5 $ minutes are constructed by averaging over each consecutive $ 12 $ seconds within each five minutes. That is, we assume that the power consumption during the $ t $-th time slot of the $ k $-th sequence is equal to the average of the $ k $-th twelve seconds within the $ t $-th time slot. As a result, we obtain $ 25 $ load sequences that have a total duration of $ 18 $ days each, and of $ 5 $-minute resolution. Similarly, we use the PV generation data available at~\cite{bank2013development} to obtain similar training data representing generation profiles from renewable energy sources. 

\begin{figure}[!t]
	\centering
	\includegraphics[width=0.45\textwidth]{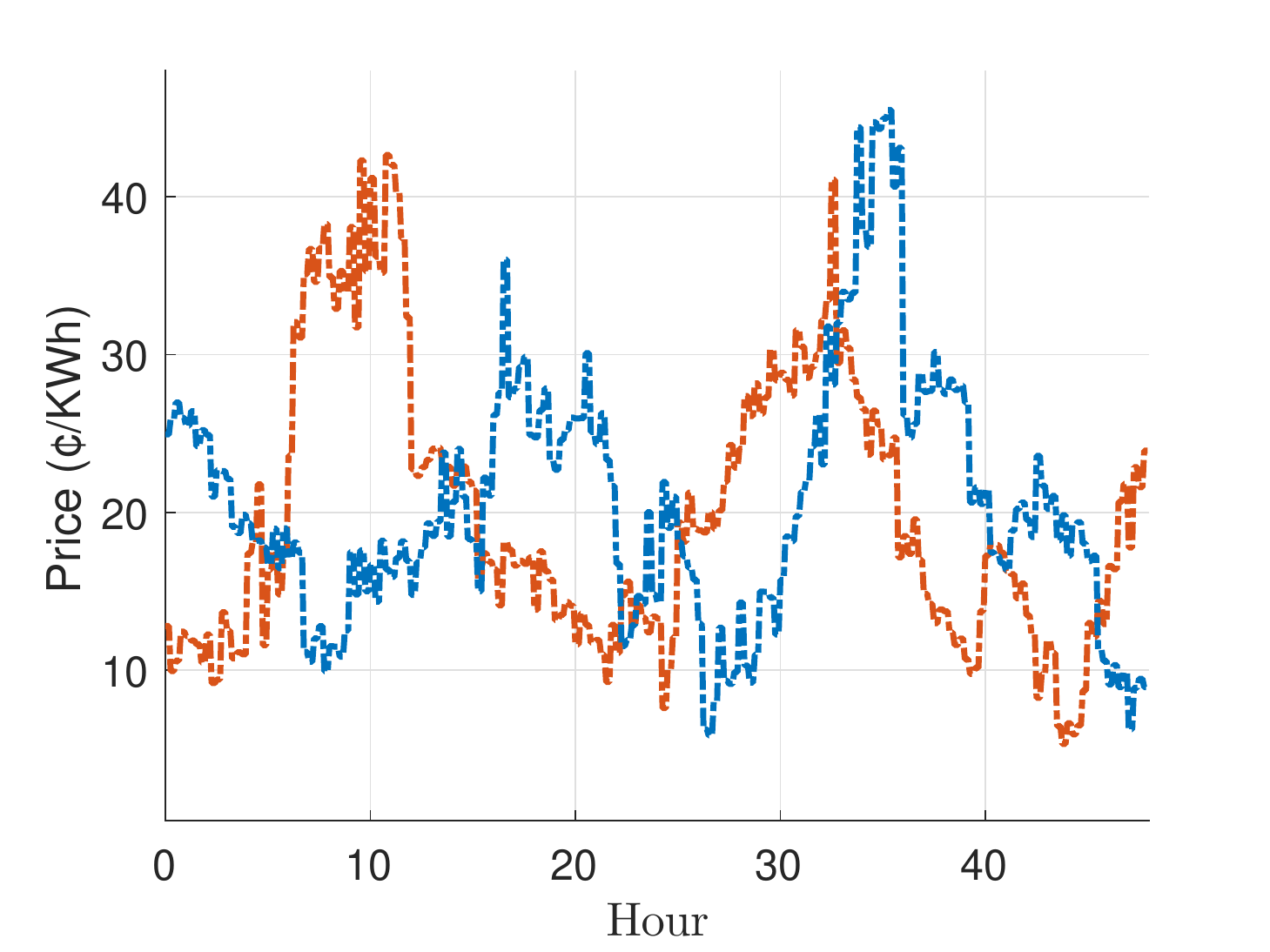}
	\caption{Two examples of energy price sequences.}
	\label{fig:Prices}
\end{figure}

For the energy prices $ p_t $, we utilized the discrete Markov processes approach~\cite{Olsson-2008} to design a Markov process that generates price sequences that has the characteristics of real prices available from the Midcontinent Independent System Operator\footnote{Available at: \url{https://www.misoenergy.org}}. The Markov chain has nine states associated with different mean values. The price of electricity at any time slot is modeled as random variable with the mean associated with the current Markov process state. In addition, we assume that the price changes happen every $ 15 $ minutes, i.e., three time slots. Fig.~\ref{fig:Prices} shows three different pricing sequences generated using the approach for a period of $ 48 $ hours.

For simplicity, we assume that the utility function of the customer at time $ t $ is the negative payment to the grid. Therefore, the utility function $ u_t(b_t + d_t - r_t, p_t) $ can be formulated as
\begin{equation}
u_t(b_t + d_t - r_t, p_t) = - p_t (b_t + d_t - r_t).
\end{equation}

Our experimentations suggest that a three-layer neural network with sigmoid activation function suffices to obtain near-optimal performance. The number of nodes in the hidden layers of the neural network is $ 128 $ and $ 32 $, respectively. We use the training approach described in Algorithm~\ref{alg:DQN} to train the neural network from the historical training data using the Matlab Neural Networks toolbox.

In order to assess the performance of the proposed approach on the testing set we use two different baseline algorithms. The first baseline assumes exact prior knowledge of the energy forecast as well as the energy prices during the optimization period. The second algorithm assumes exact knowledge of the energy forecasts and prices for a small future time interval. The detailed description of both approaches follows.

\underline{Optimal solver:} assumes that the energy prices and forecast are known accurately before the operation window. Therefore, an optimization problem can formulated to minimize the consumer's total payments~\eqref{eq:Cost} subject to the storage operational constraints~\eqref{eq:rate-con} and~\eqref{eq:Dyn_battery} at all time slots. This solver provides a lower bound on the cost that the consumer can pay using any energy storage management scheme. 

\underline{Model predictive control (MPC) solver:} assumes knowledge of the energy forecasts and prices for a short period in the future, which is called $ \tau $. Therefore, an optimization problem can be solved to find the optimal charging and discharging schedule for the period of known prices and forecasts, $ \tau $. Consequently, the consumer implements the first decision in that schedule. Then, a new optimization problem is solved incorporating new information about energy forecast and price at an additional time slot. In other words, the forecast data is available in a sliding window fashion. In addition to being very expensive computationally, this approach assumes perfect knowledge of future prices which is hard to obtain when a distribution network is operated according to a real-time pricing scheme. This approach is equivalent to the optimal solver when the value of $ \tau $ is the same as the whole optimization time.

For our simulations, we set the size of the energy storage unit to be $ E = 200 $ KWh, and the maximum charging and discharging rates are set to be $ 50  $ KW. The inefficiency coefficients $ \eta_c $ and $ \eta_d $ are set to be $ 0.9 $.

We train the proposed DQN using $ 23 $ sequences of length $ 18 $ days each, and we use the remaining $ 2 $ sequences to test the performance of the proposed approach against the two baselines. After training the DQN, the agent operates the storage unit according to the action with the highest expected reward, i.e.,  $\arg\max_{a'} Q({\bf s}, a') $. We divide testing data into eighteen intervals with $ 48 $ hours each and a resolution of $ 5 $ minutes. Then, we compare the cost of operating the storage unit using our approach, the optimal approach, and the MPC solver with $ \tau = 1 $ hour. Table~\ref{tab:res1} shows that the proposed approach achieves an average operational cost close to the optimal cost without any knowledge about the forecasted energy generation.

\begin{table}[h]
	\centering
	\begin{tabular}{| c | c|} 
		\hline
		{\bf Algorithm} &  {\bf Cost (\$) } \\ 
		\hline\hline
		{Optimal} &  366.7134   \\ 
		\hline
		{MPC} & 427.8700    \\
		\hline
		{DQN} & 429.0287 \\
		\hline
	\end{tabular}
	\caption{Average operational cost over $ 48 $ hours using different algorithms.}
	\label{tab:res1}
\end{table}

\begin{figure}[!t]
	\centering
	\begin{subfigure}[b]{0.5\textwidth}
	\centering
		\includegraphics[width=0.75\textwidth]{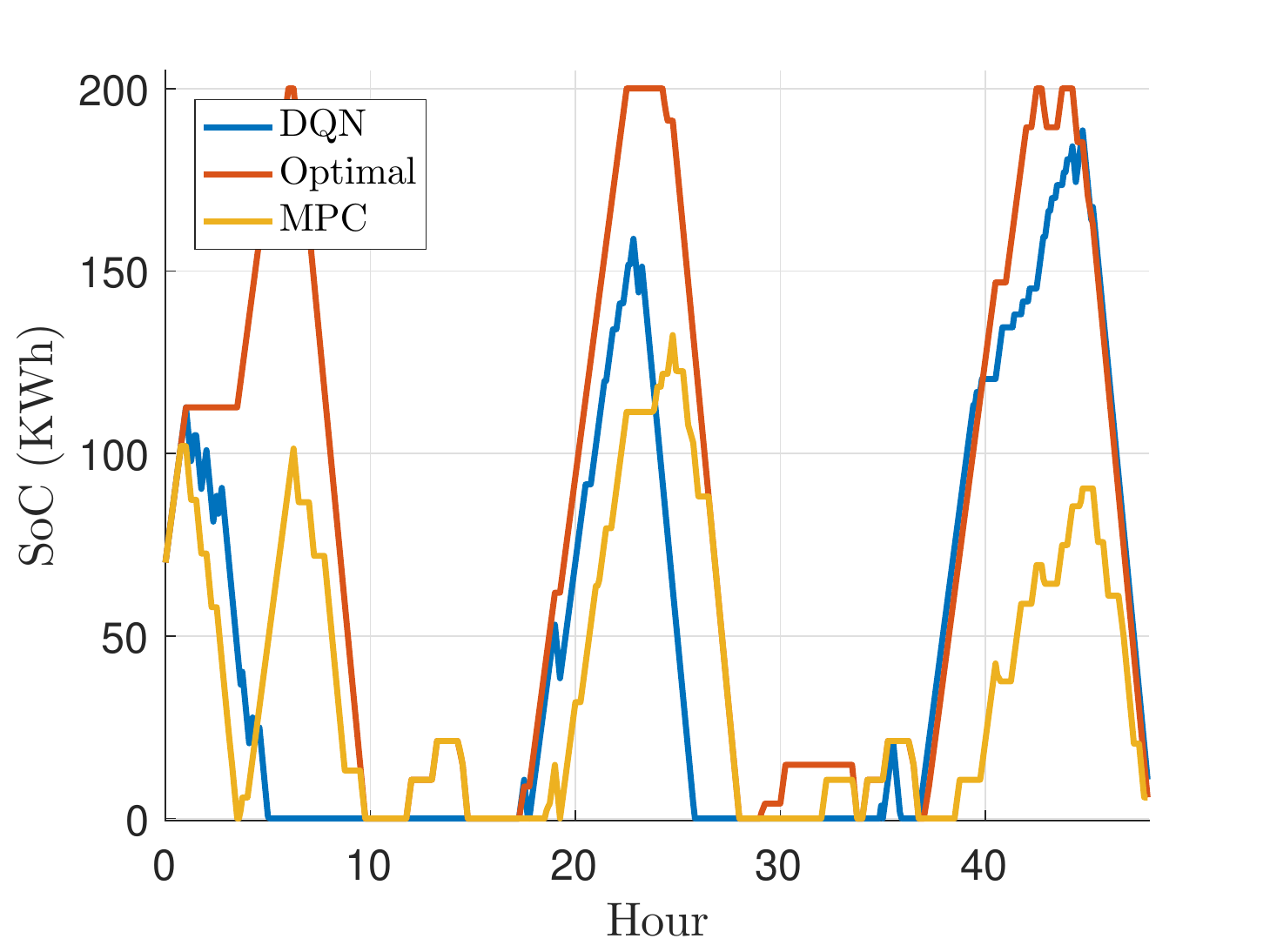}
	\end{subfigure}

	\begin{subfigure}[b]{0.5\textwidth}
	\centering
	\includegraphics[width=0.75\textwidth]{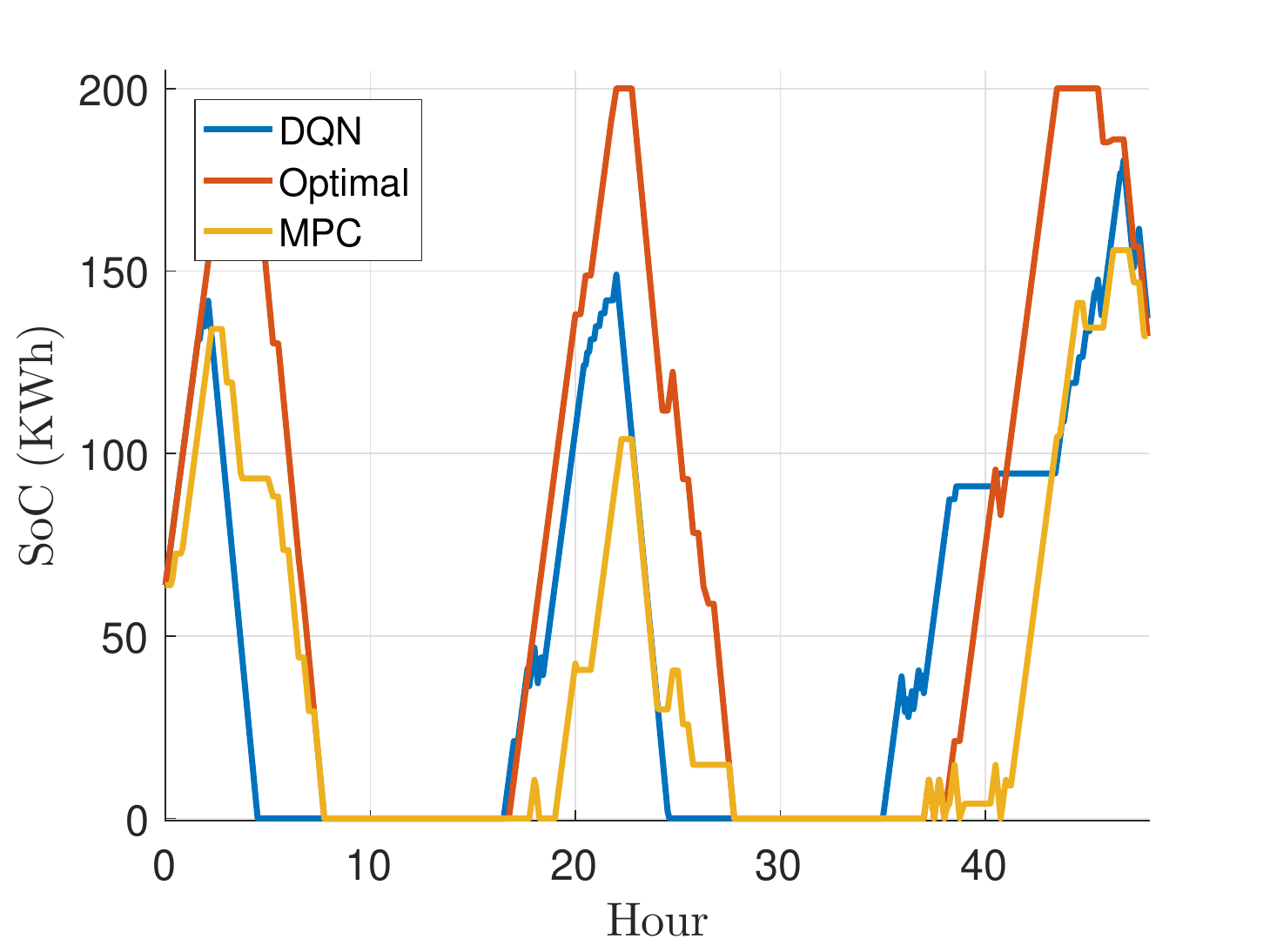}
	\end{subfigure}
	\caption{Two examples of the battery state of charge when controlled using different policies.}
	\label{fig:SoC}
\end{figure}

In Fig.~\ref{fig:SoC}, we compare the battery state of charge using our approach against the two baselines for two different scenarios. The value of $ \tau $ is set to be $ 2 $ hours for this experiment. It is clear that the battery utilization pattern using our DQN-based approach follows the pattern of the optimal approach with no knowledge of the future energy forecasts or prices. In addition, the storage utilization pattern using the MPC approach also follows the same trend but at a much higher complexity where an optimization problem needs to be solved at each time slot. 

\vspace{10pt}
\section{Conclusions}
This paper presents a novel approach for the control of energy systems (or microgrids) comprising energy storage device, renewable energy source, and inelastic load using deep Q-network. The approach learns an optimal policy by interacting with the environment. In order to approximate the action-value function, the proposed approach utilizes a neural network where the parameters are estimated from observing the rewards and state transitions associated with implemented actions. The merits of the proposed approach were demonstrated using semi-sythetic scenarios where near-optimal solutions are achieved with no distributional assumptions.

\newpage

\bibliographystyle{IEEEtran}
\bibliography{IEEEabrv,DSM}
\noindent



\end{document}